\input amstex
\documentstyle{amsppt}
\magnification=\magstep1
 \hsize 13cm \vsize 18.35cm \pageno=1
\loadbold \loadmsam
    \loadmsbm
    \UseAMSsymbols
\topmatter
\NoRunningHeads
\title Symmetry Properties of the generalized higher-order Euler polynomials
 \endtitle
\author
  Taekyun Kim \endauthor
 \keywords Euler numbers and polynomials, Symmetry
\endkeywords

\abstract
The purpose of this paper is to generalize this relation of symmetry between the power sum polynomials  and the generalized
Euler polynomials to the relation between the power sum polynomials and the generalized  higher-order  Euler polynomials.
\endabstract
\thanks  2000 AMS Subject Classification: 11B68, 11S80
\newline  The present Research has been conducted by the research
Grant of Kwangwoon University in 2009
\endthanks
\endtopmatter

\document

{\bf\centerline {\S 1. Introduction}}

 \vskip 20pt
Let $d$ be a fixed positive odd integer and let $\chi$ be the Dirichlet's character with conductor $d$.
Then the generalized Euler numbers and polynomials attached to $\chi$ are defined as
$$\aligned
&\frac{2\sum_{a=0}^{d-1}(-1)^a \chi(a)e^{at}}{e^{dt}+1}=\sum_{n=0}^{\infty}E_{n,\chi}\frac{t^n}{n!}, \\
&\frac{2\sum_{a=0}^{d-1}(-1)^a \chi(a)e^{at}}{e^{dt}+1}e^{xt}=\sum_{n=0}^{\infty}E_{n,\chi}(x)\frac{t^n}{n!}
 ,\text{ for $|t|<\frac{\pi}{d}$, (see  [1-14])}.
\endaligned\tag1$$
  For a real or complex parameter $\alpha$,we define the generalized higher-order Euler numbers and polynomials attached to $\chi$ as follows:
  $$\aligned
  &\left( \frac{2\sum_{a=0}^{d-1}(-1)^a\chi(a)e^{at}}{e^{dt}+1}\right)^{\alpha}=\sum_{n=0}^{\infty}E_{n,\chi}^{(\alpha)}\frac{t^n}{n!}, \text{ and }\\
  & \left( \frac{2\sum_{a=0}^{d-1}(-1)^a\chi(a)e^{at}}{e^{dt}+1}\right)^{\alpha}e^{xt}=\sum_{n=0}^{\infty}E_{n,\chi}^{(\alpha)}(x)\frac{t^n}{n!},
  \text{ where $|t|<\frac{\pi}{d}$.}
  \endaligned\tag2$$
  From (2), we note that
  $$E_{n,\chi}^{(\alpha)}(x)=\sum_{l=0}^n\binom{n}{l}E_{l,\chi}^{(\alpha)} x^{n-l}. \tag3$$
  By (1), we can also easily see that
  $$E_{k,\chi}(nd)+E_{k,\chi}=2T_{k, \chi}(nd-1), \text{ where $T_{k,\chi}(n)=\sum_{l=0}^n (-1)^l\chi(l)l^k $.} \tag4$$
  Recently, Kurt (see [10]) and  several authors (see [5-8], [11-13]) have studied the symmetry property for
  the Bernoulli and Euler polynomials.
  The main purpose of this paper is to prove  an identity of symmetry for the generalized higher-order Euler polynomials
  by using the method of Kurt. It turn out that the recurrence relation
  and multiplication theorem for the generalized Euler polynomials attached to $\chi$.

\vskip 20pt

{\bf\centerline {\S 2. Identities related to the generalized higher-order Euler polynomials}} \vskip 10pt
    For $w_1, w_2 \in \Bbb N$ with $w_1 \equiv 1, w_2\equiv 1$ $ (
    \mod 2)$, let us consider the following functional equation:
 $$\aligned
  I&=\frac{1}{2}\left( \frac{2\sum_{a=0}^{d-1}(-1)^a\chi(a)e^{w_1at}}{e^{dw_1t}+1}\right)^{m}e^{w_1w_2xt}\left(e^{dw_1w_2t}+1\right) \\ &\times \left( \frac{2\sum_{a=0}^{d-1}(-1)^a\chi(a)e^{w_2at}}{e^{dw_2t}+1}\right)^{m}e^{w_1w_2yt}\\
  &=\left( \frac{2\sum_{a=0}^{d-1}(-1)^a\chi(a)e^{w_1at}}{e^{dw_1t}+1}\right)^{m}e^{w_1w_2xt}\left( \frac{e^{dw_1w_2t}+1}{e^{w_2dt}+1} \right)
  \left( \sum_{b=0}^{d-1}\chi(b)(-1)^b e^{w_2bt}  \right)\\
  & \times \left(  \frac{2\sum_{a=0}^{d-1}(-1)^a\chi(a)e^{w_2at}}{e^{dw_2t}+1}\right)^{m-1}e^{w_1w_2yt}.
   \endaligned\tag5$$
It is easy to see that
$$\aligned
&\left( \frac{ e^{w_1dt}+1}{e^{dt}+1} \right)\left(   \sum_{i=0}^{d-1}\chi(i)e^{it}  (-1)^i\right)
=\left(\sum_{l=0}^{w_1-1}e^{ldt}(-1)^l \right)\left(\sum_{i=0}^{d-1}\chi(i)(-1)^ie^{it}\right)\\
&=\sum_{i=0}^{w_1d-1}\chi(i)(-1)^ie^{it}=\sum_{k=0}^{\infty}\left( \sum_{i=0}^{w_1d-1}\chi(i)(-1)^ii^k \right)\frac{t^k}{k!}
=\sum_{k=0}^{\infty}T_{k,\chi}(w_1d-1)\frac{t^k}{k!}.\\
\endaligned\tag6$$
By (5) and (6), we have
$$\aligned
&I=\left(\sum_{i=0}^{\infty}E_{i,\chi}^{(m)}(w_2x)\frac{w_1^it^i}{i!}\right)\left(\sum_{l=0}^{\infty}T_{l,\chi}(w_1d-1)\frac{w_2^l t^l}{l!}\right)
\left(\sum_{k=0}^{\infty}E_{k, \chi}^{(m-1)}(w_1y)\frac{w_2^kt^k}{k!}\right)\\
&=\sum_{n=0}^{\infty}\left(\sum_{j=0}^n \binom{n}{j}w_2^jw_1^{n-j}E_{n-j,\chi}^{(m)}(w_2x)\sum_{k=0}^j T_{k,\chi}(w_1d-1)\binom{j}{k}E_{j-k,\chi}^{(m-1)}(w_1y)\right)\frac{t^n}{n!}.
\endaligned\tag7$$
By the symmetry of $I$ in $w_1$ and $w_2$, we also see that
$$I=\sum_{n=0}^{\infty}\left(\sum_{j=0}^n\binom{n}{j}w_1^jw_2^{n-j}E_{n-j,\chi}^{(m)}(w_1x)\sum_{k=0}^j\binom{j}{k}T_{k,\chi}(w_2d-1)
E_{j-k,\chi}^{(m-1)}(w_2y)\right)\frac{t^n}{n!}. \tag8$$

By comparing the coefficients on the both sides of (7) and (8), we obtain the following theorem.

\proclaim{ Theorem 1} Let $\chi$ be the Dirichlet's character with
an odd conductor $d\in\Bbb N$. For  $w_1, w_2, d\in \Bbb N$ with $d
\equiv 1 (\mod 2)$, $w_1\equiv 1(\mod 2)$, $w_2\equiv 1 (\mod 2)$,
we have
$$\aligned
&\sum_{j=0}^n \binom{n}{j}w_2^jw_1^{n-j}E_{n-j, \chi}^{(m)}(w_2 x) \sum_{k=0}^j\binom{j}{k}T_{k, \chi}(w_1d-1)E_{j-k,\chi}^{(m-1)}(w_1y)\\
& =\sum_{j=0}^n\binom{n}{j}w_1^jw_2^{n-j}E_{n-j,\chi}^{(m)}(w_1x)\sum_{k=0}^j\binom{j}{k}T_{k,\chi}(w_2d-1)
E_{j-k,\chi}^{(m-1)}(w_2y).\endaligned\tag9$$
\endproclaim

Let $y=0$ and $m=1$ in (9). Then we have  the following corollary.

\proclaim{ Corollary 2}
For $ n\in\Bbb Z_+,$ we have
$$\aligned
&\sum_{j=0}^n\binom{n}{j}w_2^jw_1^{n-j}E_{n-j,\chi}(w_2x)T_{j,\chi}(w_1d-1)\\
&=\sum_{j=0}^n\binom{n}{j}w_1^jw_2^{n-j}E_{n-j,\chi}(w_1x)T_{j,\chi}(w_2d-1).
\endaligned\tag10$$
\endproclaim

From (5) and (6), we can also derive  the following Eq.(11).
$$\aligned
&I=\left( \frac{2\sum_{a=0}^{d-1}\chi(a)(-1)^a e^{w_1at}}{e^{dw_1t}+1} \right)^m e^{w_1w_2xt}
\left( \sum_{i=0}^{w_1d-1}\chi(i)(-1)^ie^{iw_2t} \right)\\
&\times\left(\frac{2\sum_{b=0}^{d-1}\chi(b)(-1)^be^{w_2bt}}{e^{dw_2t}+1}\right)^{m-1} e^{w_1w_2yt}\\
&=\left( \sum_{i=0}^{w_1d-1}\chi(i)(-1)^i\sum_{k=0}^{\infty}E_{k,\chi}^{(m)}(w_2x+\frac{w_2}{w_1}i)\frac{w_1^kt^k}{k!}\right)
\left(  \sum_{l=0}^{\infty}E_{l,\chi}^{(m-1)}(w_1y)\frac{w_2^lt^l}{l!} \right)\\
&=\sum_{n=0}^{\infty}\left(\sum_{k=0}^n \binom{n}{k}w_1^kw_2^{n-k}E_{n-k,\chi}^{(m-1)}(w_1y)\sum_{i=0}^{w_1d-1}\chi(i)(-1)^iE_{k,\chi}^{(m)}(w_2x+\frac{w_2}{w_1}i)\right)
\frac{t^n}{n!}.
\endaligned\tag11$$
By the symmetry of $I$ in $w_1, w_2$, we also see that
$$I=\sum_{n=0}^{\infty}\left(\sum_{k=0}^n \binom{n}{k}w_2^kw_1^{n-k}E_{n-k,\chi}^{(m-1)}(w_2y)\sum_{i=0}^{w_2d-1}\chi(i)(-1)^iE_{k,\chi}^{(m)}
(w_1x+\frac{w_1}{w_2}i)\right)\frac{t^n}{n!}.
\tag12$$
By comparing the coefficients on the both sides of (11) and (12), we obtain the following theorem.

\proclaim{ Theorem 3}
Let $w_1, w_2$ be the odd natural numbers. For $n\in\Bbb Z_{+}, m\in\Bbb N$, we have
$$\aligned
&\sum_{k=0}^n \binom{n}{k}w_1^kw_2^{n-k}E_{n-k,\chi}^{(m-1)}(w_1y)\sum_{i=0}^{w_1d-1}\chi(i)(-1)^i
E_{k,\chi}^{(m)}(w_2x+\frac{w_2}{w_1}i)\\
&=\sum_{k=0}^n \binom{n}{k}w_2^kw_1^{n-k}E_{n-k,\chi}^{(m-1)}(w_2y)\sum_{i=0}^{w_2d-1}\chi(i)(-1)^iE_{k,\chi}^{(m)}
(w_1x+\frac{w_1}{w_2}i)\\
\endaligned\tag13$$
\endproclaim
  Let $m=1$ in (13). Then we have
   $$w_1^n\sum_{i=0}^{w_1d-1}\chi(i)(-1)^i
 E_{n,\chi}(w_2x+\frac{w_2}{w_1}i)=w_2^n\sum_{i=0}^{w_2d-1}\chi(i)(-1)^iE_{n,\chi}(w_1x+\frac{w_1}{w_2}i)$$

\vskip 20pt

\quad Taekyun Kim

\quad Division of General Education-Mathematics, Kwangwoon
University, Seoul

\quad 139-701, S. Korea
 e-mail:\text{ tkkim$\@$kw.ac.kr}

\enddocument